\documentclass{article}

\usepackage[utf8]{inputenc}

\usepackage{amsfonts}
\usepackage{amssymb}
\usepackage{amsmath}
\usepackage{amssymb}
\usepackage{float}
\usepackage{amsthm}
\usepackage[russian]{babel}
\usepackage[pdftex]{graphicx}

\usepackage{fullpage}
\usepackage{wrapfig}

\newtheorem{theorem}{Теорема}
\newtheorem{question}{Вопрос}
\newtheorem{prop}{Утверждение}
\newtheorem{corr}{Следствие}
\newtheorem{lemma}{Лемма}
\newtheorem{definition}{Определение}

\def\Int{\operatorname{Int}}

\def\Diam{\operatorname{Diam}}

\begin {document}


\title{{\small УДК 514.17+519.17+ 515.124.4} О хроматическом числе плоскости.}
\author{А.~Я.~Канель-Белов, В.~А.~Воронов, Д.~Д.~Черкашин}

\maketitle


\begin {abstract}

Данная статья посвящена естественному обобщению задачи о хроматическом числе плоскости.
Рассматриваются хроматические числа пространств вида $\mathbb{R}^n \times [0, \varepsilon]^k$ с~запрещенным евклидовым расстоянием $1$.

Показано, что  $5 \leq\chi(\mathbb{R}^2\times [0,\varepsilon])\leq 7$ и ${6\leq\chi(\mathbb{R}^2 \times [0,\varepsilon]^2)\leq 7}$
при достаточно малом $\varepsilon > 0$.

Также в статье рассматриваются естественным образом возникающие дальнейшие вопросы.

\end {abstract}

\medskip
{\bf Ключевые слова:}\ Хроматическое число плоскости, хроматические числа пространств.





\begin{abstract}
This paper is devoted to a natural generalization of the chromatic number of the plane problem.
We consider chromatic number of the spaces $\mathbb{R}^n \times [0,\varepsilon]^k$.

We prove that  $5 \leq\chi(\mathbb{R}^2\times [0,\varepsilon])\leq 7$ and ${6\leq \chi(\mathbb{R}^2\times [0,\varepsilon]^2) \leq 7}$.

We also observe natural questions arising from these considerations.
\end{abstract}

\medskip


\section {Введение}

Рассмотрим граф, вершинами которого являются точки плоскости, а ребра соединяют пары точек на расстоянии $1$.
Э.~Нельсон поставил задачу о нахождении хроматического числа этого графа (обозначим его за $\chi (\mathbb{R}^2)$). Затем задача была популяризована М.~Гарднером, П.~Эрдешом, Г.~Хадвигером и А.~Сойфером. В русскоязычной литературе устоялось название {\it проблема Нельсона--Хадвигера}. Мы благодарим А.~Сойфера за указание исторических неточностей.

Хорошо известна следующая теорема:

\begin{theorem}            \label{Th1}
$4 \leq \chi(\mathbb{R}^2) \leq 7$.
\end{theorem}

Получить эти оценки сравнительно несложно.
К сожалению, более чем за 65 лет не появилось никаких аргументов, позволяющих их улучшить. С другой стороны, попытки как-то подойти к этой проблеме породили
огромное количество интересных задач и содержательных результатов, им посвящена обширная литература (чуть более подробная  библиография --- см. раздел \ref{Obzor}).

Естественным ослаблением служит требование нахождения точек на расстоянии, сколь угодно близком к единице. В этом случае найдутся две точки на почти единичном расстоянии для любой раскраски плоскости в $5$ цветов, см. например~\cite{currie2015chromatic} (более слабые утверждения есть в~\cite{exoo2005varepsilon} и~\cite{grytczuk2016fractional});
это утверждение является очевидным следствием Теоремы 9, которая доказывается далее.
Более сильный результат состоит в том, что если одноцветные множества измеримы, то при четырехцветной раскраске также найдутся точки на единичном расстоянии (см.~\cite{F}).

Мы рассматриваем ``почти плоский'' случай или случай размерности ``$2+\varepsilon$''. Наш результат состоит в том, что если слойка между двумя плоскостями в трехмерном пространстве {\it произвольным образом} раскрашена в $4$ цвета, то найдутся две точки на единичном расстоянии (см. теорему \ref{Th6}).

Кроме того, мы показываем что при наличии двух инфинитезимальных измерений, это верно и для пяти цветов. Иными словами, если прямое произведение плоскости на сколь угодно малый квадрат раскрашено в $5$ цветов, то найдутся две одноцветные точки на единичном расстоянии  (см. теорему \ref{Th7}).

Случай прямого произведения выглядит существенно интереснее и сложнее случая расстояний из интервала $(1-\varepsilon,\ 1+\varepsilon)$. Например, нам удается найти точки на единичном расстоянии в слойке только для четырех цветов. В этой связи возникает естественный

\begin{question}
Пусть прямое произведение плоскости на отрезок раскрашено в $5$ цветов. Верно ли, что найдутся две одноцветных точки на единичном расстоянии?
\end{question}

Другим направлением исследований явилось следующее.

\begin{question}
Рассмотрим раскраску плоскости в $n$ цветов. Верно ли, что в одном из цветов укладываются все расстояния?
\end{question}

Иными словами, может ли так быть, чтобы в первом цвете не наблюдалось расстояния $d_1$, во втором~--- $d_2$, в третьем~--- $d_3$? Когда $d_i=1$ мы получаем задачу о хроматическом числе плоскости. Интересно, что для трех цветов данная задача оказалась весьма нетривиальной, а для $6$ были построены различные контрпримеры.

И в этой связи возникают аналогичные вопросы для слоек и почти расстояний.


\section{О хроматических числах пространств}  \label{Obzor}

Соответствующая тематика активно развивается, получено много интересных результатов и оценок в~различных ситуациях (например, см.~\cite{Soi2}). Мы приведем лишь краткий перечень основных достигнутых результатов. Хроматические числа пространств активно исследовались, например, в школе А.~М.~Райгородского. Более подробные сведения о проблеме Нельсона--Хадвигера и смежных задачах можно почерпнуть из следующих обзоров: П.~К.~Агарвал и Я.~Пах~\cite{AP}, П.~Брасс, В.~Мозер, Я.~Пах \cite{BMP}, М.~Бенда, М.~Перлес~\cite{BP}, К.~Б.~Чилакамарри~\cite{Ch}, В.~Кли и С.~Вэгон~\cite{KW}, А.~М.~Райгородский~\cite{Rai1},~\cite{Rai2},~\cite{Rai8},~\cite{Rai9},~\cite{Rai10},  А.~Сойфер~\cite{Soi2,Soi3} и Л.~А.~Секеи~\cite{Szek}.

\subsection{О хроматическом числе плоскости}
Начнем с ослаблений. Если каждое одноцветное множество разбивается на связные области, ограниченные жордановыми кривыми, то необходимо не менее 6 цветов, что было доказано Д.~Р.~Вудаллом еще в 1973-ем году~\cite{W}. К.~Дж.~Фалконер в 1981-ом году показал, что если потребовать, чтобы множества точек, раскрашенных в один и тот же цвет, были измеримы по Лебегу, тогда для правильной раскраски плоскости требуется хотя бы $5$ цветов~\cite{F}.
Разумеется, раскраска плоскости в $7$ цветов обеспечивает оценку сверху и для ослабленных формулировок
(см. Рис. 1).

\begin{figure}[ht]
  \centering
  \includegraphics[width=6cm]{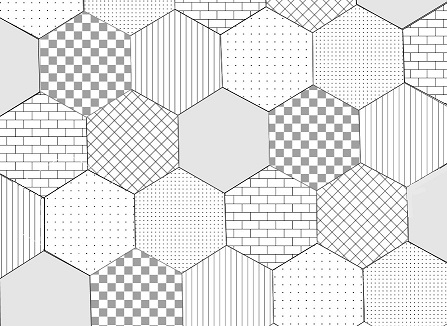}\\
  \caption{Раскраска плоскости в $7$ цветов. Стороны правильных $6$-угольников имеют длину $\frac{1}{\sqrt{7}}$.}
\end{figure}

Одна из главных трудностей заключается в том, что ответ \textit{может зависеть} от те\-о\-ре\-ти\-ко-мно\-жес\-твен\-ной аксиоматики, как показали в 2003-ем году C.~Шелах и А.~Сойфер~\cite{ShS}.
Если мы предполагаем аксиому выбора, то по теореме Эрд\-е\-ша-де Брейна~\cite{EB} хроматическое число бесконечного графа реализуется на конечном подграфе. Однако компьютерный перебор не находит подграфов с хроматическим числом хотя бы 5, что позволяет предположить, что хроматическое число в стандартной аксиоматике равняется четырем.
Если же отказаться от аксиомы выбора, но дополнить стандартную аксиоматику Цермело-Френкеля аксиомой зависимого выбора и дополнительно потребовать измеримость всех
подмножеств по Лебегу, то применимо доказательство Фалконера, и хроматическое число лежит между пятью и семью.

Также можно рассматривать случай ограниченного подмножества плоскости, что было сделано, например, в~\cite{Kr}.

\subsection{Задача для произвольных метрических пространств}

Перейдем к различным обобщениям. Для произвольного метрического пространства $(X, d)$ и числа $a > 0$ определим граф $G_a(X,d)$ следующим образом:
множество его вершин совпадает с точками пространства, множество ребер образуют пары точек, лежащие на расстоянии $a$. Нас по-прежнему интересует хроматическое число этого графа $\chi(X,d,a)$.
Наиболее часто в качестве $(X,d)$ рассматривают $\mathbb{R}^n$ и $\mathbb{Q}^n$ c евклидовой метрикой. Мы ограничимся этими случаями при $a=1$.
Стоит отметить, что в~вещественном случае все графы $G_a$, очевидно, изоморфны.

\subsubsection {Хроматические числа вещественных пространств}

Для прямой ответ очевиден: $\chi(\mathbb{R}) = 2$. При $n =2$ мы получаем в точности исходную формулировку вопроса о хроматическом числе плоскости.

При $n=3$ задача представляется еще более сложной, чем классическая проблема Нель\-со\-на--Хад\-ви\-ге\-ра,~--- последние оценки с обеих сторон получены в текущем столетии.

\begin{theorem}
$$6 \leq \chi(\mathbb{R}^3) \leq 15.$$
\end{theorem}
Нижняя оценка принадлежит О.~Нечуштану~\cite{Nech}, верхняя~--- Д.~Кулсону~\cite{Coul}.

В асимптотике же выполняются следующие оценки:

\begin{theorem}
 $$\left(1.239\ldots+o(1)\right)^n \le \chi\left({\mathbb R}^n\right) \le (3+o(1))^n.$$
 \end{theorem}

Нижняя оценка принадлежит А.~М.~Райгородскому~\cite{Rai7}, верхняя~--- Д.~Г.~Ларману и К.~А.~Роджерсу~\cite{LR}. Стоит отметить, что асимптотические нижние оценки в этой и смежных задачах были
получены интересным самим по себе \textit{линейно-алгебраическим методом}, более того, с их помощью Дж.~Кан и Дж.~Калаи в 1993-ем году построили контрпример к~гипотезе Борсука~\cite{KK}, стоявшей к~тому моменту более пятидесяти лет.
Более подробно про гипотезу Борсука и линейно-алгебраический метод написано в~\cite{Rai8}.

\subsubsection {Хроматические числа рациональных пространств}

Одномерный случай, как и для вещественных чисел, тривиален:  $\chi(\mathbb{Q})=2$.

Некоторое удивление может вызвать тот факт, что точное значение хроматического числа $\mathbb{Q}^n$ известно не только в размерности $2$, но и в размерностях $3$ и $4$
(см. Д.~Р.~Вудалл~\cite{W}, П.~Д.~Джонсон~\cite{J} и М.~Бенда-М.~Перлес~\cite{BP}).

\begin{theorem}       \label{Th4}
 $\chi(\mathbb{Q}^2) = \chi(\mathbb{Q}^3) = 2$, $\chi(\mathbb{Q}^4) = 4$.
\end{theorem}

Лучшая асимптотическая нижняя оценка принадлежит Е.~И.~Пономаренко и А.~М.~Райгородскому~\cite{PR1},~\cite{PR2}, верхняя~--- Д.~Г.~Ларману и К.~А.~Роджерсу~\cite{LR}.

\begin{theorem}            \label{Th5}
$$\left(1.199\ldots + o(1)\right)^n \le \chi\left({\mathbb Q}^n\right) \le \left(3+o(1)\right)^n.$$
\end{theorem}

Отметим также, что в статье~\cite{Ax} рассматривался смешанный случай пространства $\mathbb{R} \times \mathbb{Q}$.

\subsection{Полихроматическое число}

В книге~\cite{HD} Г.~Хадвигер и Х.~Дебрюннер с подачи П.~Эрдеша сформулировали естественный вопрос об отыскании \textit{полихроматического числа плоскости}, заключающийся в нахождении
наименьшего числа цветов, необходимых для такой раскраски плоскости, что для любого цвета $i$ найдется расстояние $d$, такое что цвет $i$ не содержит пару точек на расстоянии $d$.
Мы используем обозначение $\chi_p$ для вышеупомянутой величины, предложенное Сойфером в статье~\cite{Soi}.
Наилучшие оценки были получены Райским и Стечкиным в работе~\cite{Rai} (пример Стечкина опубликован с его разрешения в той же статье).

\begin{theorem}
$$4 \leq \chi_p \leq 6.$$
\end{theorem}

\noindent Другие доказательства тех же оценок можно найти в статье Д.~Р.~Вудалла~\cite{W}.

\section{Основные результаты}

Нас интересуют хроматические числа пространств вида $\mathbb{K}^n \times [0,\varepsilon]^k$, где $\mathbb{K}\in\{\mathbb{R},\mathbb{Q}\}$, $n,k \geq 1$ с евклидовой метрикой. Такие метрические пространства будем называть ``слойками'', а в случае $n = k = 1$ ``полосами''.
В данной работе основное внимание уделяется случаю $n=2$.

\subsection {Хроматические числа одномерных слоек}

Начнем с простого наблюдения. Нижняя оценка есть в~\cite{bauslaugh1998tearing} и~\cite{Ch}, но мы приведем ее для полноты изложения.

\begin{prop} \label{Utv1}
 Пусть $0 < h\leq \sqrt{\frac{3}{4k}}$. Тогда $$\chi(\mathbb{R} \times [0,h]^k) = 3.$$
 Пусть $\sqrt{\frac{3}{4k}} < h \leq \sqrt{\frac{8}{9k}}$. Тогда $$\chi(\mathbb{R}\times [0, h]^k) = 4.$$
\end{prop}

\noindent{\bf Верхняя оценка.}\ Пусть $0 < h \leq\sqrt{\frac{3}{4k}}$. Раскрасим $\mathbb{R}$ в $3$ цвета, чередуя одноцветные полуинтервалы длины $1/2$ (цветов $1, 2, 3, 1, 2, 3$ и т.д.). Затем каждой точке $\mathbb{R}\times [0,h]^k$ присвоим тот же цвет, который имеет проекция из прямого произведения  на действительную прямую. Тогда диаметр одноцветного параллелепипеда $[0; \frac{1}{2}] \times [0,h]^k$ не превосходит $1$, причем в случае равенства концы диаметра раскрашены по-разному.

Аналогично строится раскраска в $4$ цвета при ${\sqrt{\frac{3}{4k}} < h\leq\sqrt{\frac{8}{9k}}}$~--- в этом случае полуинтервалы имеют длину $1/3$.

\vskip 0.2cm

\noindent{\bf Нижняя оценка.}\ 
Проиллюстрируем на тривиальном примере схему доказательства, которая будет использована далее для размерностей 3 и 4.

Предположим, что полоса $\mathbb{R} \times [0,\varepsilon]$ правильно раскрашена в несколько цветов. Пусть $l \in \mathbb{N}$ таково, что $1/l = \delta \leq \varepsilon^2$.  На границе слойки $\mathbb{R} \times \{0\}$ выберем раскрашенные по-разному точки $u=(x, 0)$, $v = (x+\delta,0)$, расстояние между которыми равно $\delta$. Такой выбор возможен в силу того, что точки $(0, 0)$, $(\delta,0)$, ... $(1,0)$ не могут быть все одного цвета. Обозначим $w = (x+\delta/2, \varepsilon)$. В одной из пар точек $u, w$ или $v, w$ встречаются два цвета. Пусть это пара $u, w$. Тогда для раскраски точки $\xi$, которая находится на расстоянии 1 от $u$, $w$ и лежит внутри полосы, потребуется еще один цвет. \qed

\begin{figure}[ht]
  \centering
  \includegraphics[width=6cm]{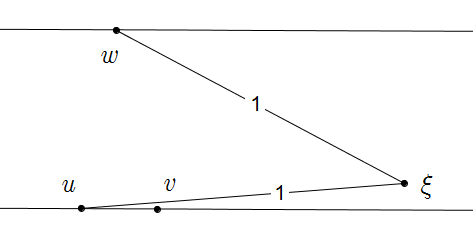}\\
  \caption{Нижняя оценка для $\chi(\mathbb{R} \times [0, \varepsilon]$).}
\end{figure}

Множество $\mathbb{R} \times [0, h]^k$ при $\sqrt{\frac{3}{4k}} < h \leq \sqrt{\frac{8}{9k}}$ содержит полосу $\mathbb{R} \times [0, h_1]$, $h_1 >\sqrt{3}/2$, которая является произведением $\mathbb{R}$ на большую диагональ $k$-мерного гиперкуба с ребром $h$.
В такую полосу вкладывается дистанционный граф, изображенный на рис. 2,  причем $d(y,z)$ может принимать любое значение из отрезка
$[0,3-2\sqrt{3-h_1^2}]$
Выберем такую реализацию этого графа, чтобы точки $x,y,z$ лежали на  границе полосы, и $d(y,z)=1/m$; $m \in \mathbb{N}$.  Копируя конструкцию $m$ раз, строим дистанционный граф с хроматическим числом $4$. \qed

\vskip 0.2cm

\textbf{Замечание.} Число вершин критического графа стремится к бесконечности, когда $h$ приближается к значению, в котором хроматическое число разрывно ($h=0$ и $h=\sqrt{\frac{3}{4k}}$), но граф может быть размещен в области, диаметр которой не зависит от $h$.

\vskip 0.2cm

Очевидно, функция $\xi_{n,k}(h) = \chi (\mathbb{R}^n\times [0,h]^k)$, определенная при $h \geq 0$, не убывает. При любых фиксированных $n$, $k$ число значений $\xi_{n,k}(h)$ конечно, поскольку $\chi(\mathbb{R}^n)\leq\xi_{n,k}(h) \leq \chi(\mathbb{R}^{n+k})$, а следовательно, конечно и число точек разрыва. По-видимому, при $n>1$ ни одна точка разрыва $\xi_{n,k}(h)$ не может быть найдена  без улучшения известных оценок $\chi(\mathbb{R}^n)$.
\vskip 0.2cm
\begin{figure}[ht]
  \centering
  \includegraphics[width=6cm]{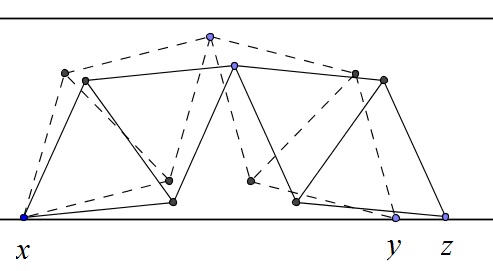}\\
  \caption{Цепочка $\theta$-графов внутри полосы.}
\end{figure}

Укажем более широкий класс множеств, для которых нижняя оценка из Утверждения \ref{Utv1} остается в силе.


\begin{prop}              \label{Utv2}
Пусть $\varepsilon$~--- произвольное положительное число; $Q$~--- $\varepsilon$-окрестность некой кривой $\xi$ диаметра хотя бы $2$.
Тогда $\chi(Q) \geq 3$.
\end{prop}


\subsection {Хроматические числа двумерных слоек}

Перейдем к промежуточному случаю между плоскостью и пространством~---
$\mathbb{R}^2 \times [0,\varepsilon]$ (слойка высоты $\varepsilon$).  Несмотря на то, что это множество по-прежнему допускает правильную раскраску в $7$ цветов, нижняя оценка менее тривиальна, чем для плоскости:

\begin{theorem}                \label{Th6}
Пусть $\varepsilon$~--- положительное число, меньшее $\sqrt{3/7}$. Тогда
$$5\leq\chi(\mathbb{R}^2\times [0,\varepsilon])\leq 7.$$
\end{theorem}

В отличие от случая полосы ($n=1$) здесь не удается доказать даже то, что функция $\chi (\mathbb{R}^2 \times [0, \varepsilon])$ разрывна в точке $\varepsilon=0$.
Рассмотрим теперь ``раздутие''\ плоскости в пространстве большей размерности.
Благодаря тому, что раскраска плоскости в $7$ цветов не содержит расстояний, принадлежащих некоторому интервалу, верхняя оценка сохраняется при увеличении размерности.

\begin{theorem}  \label{Th6prim}
Пусть $k$~--- целое число, $\varepsilon < \varepsilon_0(k)$~--- положительное число,  Тогда
$$\chi(\mathbb{R}^2\times [0,\varepsilon]^k)\leq 7.$$
\end{theorem}

Нижнюю оценку удается улучшить уже при $k=2$.

\begin{theorem}  \label{Th7}
 Пусть $\varepsilon$~--- произвольное положительное число,  Тогда $$\chi(\mathbb{R}^2\times [0,\varepsilon]^2) \geq 6.$$
\end{theorem}

Отметим, что для получения приведенных нижних оценок, как и при $n=1$, достаточно рассмотреть раскраску ограниченной области, диаметр которой не зависит от $\varepsilon$.

В доказательстве Теоремы \ref{Th5} используется следующая небезынтересная сама по себе

\begin{lemma}   \label{L1}
Предположим, что евклидова плоскость правильно покрашена в $k$ цветов. Тогда для любого $\varepsilon>0$ найдется круг радиуса  $\varepsilon$, содержащий точки как минимум трех различных цветов.
\end{lemma}

\begin{corr}  \label{Cor1}
Предположим, что евклидова плоскость правильно покрашена в $k$ цветов. Тогда для любого $\varepsilon>0$ найдется окружность радиуса меньше чем $\varepsilon$, содержащая точки как минимум трех различных цветов.
\end{corr}

Верно и более сильное утверждение, нежели Лемма~\ref{L1}:

\begin{theorem}   \label{Th8}
Пусть пространство $\mathbb{R}^n$ правильно покрашено в $m$ цветов. Иными словами, обозначив множество точек $i$-того цвета через $C_i$, имеем: $$\bigcup_{i=1}^{m} C_i = \mathbb{R}^n,$$ и ни одно из множеств $C_i$ не содержит двух точек, находящихся на расстоянии $1$. Тогда найдется $n+1$ множество из этого семейства, пересечение замыканий которых непусто.
\end{theorem}

Утверждение теоремы очевидно в том случае, если компоненты связности замыканий ${C}_i$ представляют собой многогранники,  но оно справедливо и для произвольного покрытия с одним запрещенным расстоянием.



\subsection{Хроматические числа рациональных пространств}

В рациональном случае справедлива следующая теорема:

\begin{theorem}         \label{Th9}
 Для достаточно малого положительного $\varepsilon$ выполняется $$\chi(\mathbb{Q}\times [0,\varepsilon]_\mathbb{Q}^3)=3. $$
\end{theorem}

Очевидно, нельзя заменить в условии теоремы $[0,\varepsilon]_\mathbb{Q}^3$ на $[0,\varepsilon]_\mathbb{Q}^2$ так как $\chi(\mathbb{Q}^3)=2$.



\section {Доказательства}

Начнем с конструкции, неоднократно используемой в доказательствах.

\begin{definition}
Пусть $\omega_r$~--- окружность радиуса $r$. Назовем число $r > 0$ {\em запрещенным радиусом}, если $G_1(\omega_r)$  содержит нечетный цикл.
\end{definition}


\begin{prop} \label{Utv3}
Запрещенные радиусы плотны в $[1/2, \ \infty).$
\end{prop}

\begin{proof}
В самом деле, для числа $q \in \mathbb{Q} \cap (0,\frac{1}{2})$, представимого в виде $q = \frac{l}{2k+1}$; $k,l \in \mathbb{N}$ можно построить запрещенный радиус $$r = \frac{1}{2 \sin {\pi q}}.$$
\end{proof}

\subsection {Доказательство нижней оценки в Теореме \ref{Th6}}

Будем обозначать стандартную метрику на плоскости за $\rho$, а сферу полной размерности с радиусом $r$ и центром в точке $u$ за $S(u; r)$.
Предположим, что слойка $\mathbb{R}^2 \times [0, \varepsilon]$ правильно покрашена в несколько цветов. Пусть $0<\delta < \varepsilon^2$. Выберем на границе слойки $\mathbb{R}^2 \times \{0\}$ такие по-разному раскрашенные точки $u,v$, что $\rho(u, v) =  \delta$.

\begin{figure}[ht]
  \centering
  \includegraphics[width=12cm]{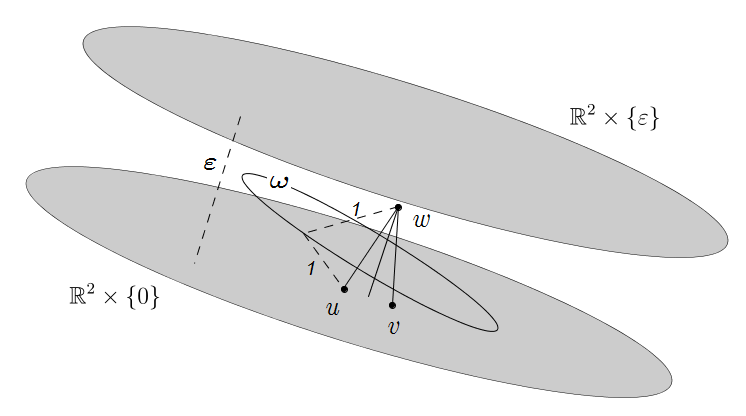}\\
  \caption{Окружность запрещенного радиуса внутри слойки $\mathbb{R}^2 \times [0, \varepsilon]$}
\end{figure}

Пусть $\varepsilon_1$ удовлетворяет условию $$\sqrt{ \delta} \leq \varepsilon_1 < \varepsilon$$
и, кроме того, $r = \sqrt{1-\varepsilon^2_1/4}$~--- запрещенный радиус. Построим равнобедренный треугольник $u v w$, высота которого $w w_1$ перпендикулярна границе слойки, а боковые стороны имеют длину $\varepsilon_1$. Поскольку $u$ и $v$ покрашены в разные цвета, хотя бы в одной из пар $u, w$ и $v, w$ встречаются два цвета. Без ограничения общности предположим, что точки $u, w$ раскрашены в цвета 1 и 2 соответственно. Тогда окружность
$$\omega = S(u; 1) \cap S(w; 1)$$
лежит в слойке, не содержит точек цветов 1 и 2, и имеет запрещенный радиус $r$, то есть для ее правильной раскраски требуется еще 3 цвета. \qed

\subsection {Доказательство верхней оценки в Теоремах~\ref{Th6} и~\ref{Th6prim}}

Рассмотрим стандартную правильную раскраску плоскости в 7 цветов (см. Рис. 1). В ней нет двух точек, покрашенных в один и тот же цвет, и лежащих на расстоянии между $2/\sqrt{7}$ и $1$.
Покрасим слойку следующим образом: каждый куб ${(x,y)} \times [0, \varepsilon]^k$ получит тот же цвет, что и точка $(x,y)$ в раскраске плоскости. Такая раскраска слойки будет правильной при условии $(2/\sqrt{7})^2 + k\varepsilon^2 < 1$, которое эквивалентно неравенствам в формулировках теорем. \qed

\subsection {Доказательство Утверждения~\ref{Utv2}}
Не умаляя общности можно считать, что $\varepsilon < 1$.
Предположим противное: существует раскраска множества $Q$ в два цвета. Обозначим за $G(Q)$ соответствующий граф и найдем в нем нечетный цикл.

\vskip+0.2cm


\vskip+0.2cm

Рассмотрим некоторую точку $u \in \xi$. Поскольку $\Diam{\xi} \geq 2$, пересечение $S(u;1)$ и $\xi$ непусто. Пусть $v \in S(u;1) \cap \xi$;
$\lVert u-v_1 \rVert = 1$; $\lVert v_i-v_{i+1} \rVert = 1$; $i = 1,2,3$.
Если углы между соседними единичными звеньями ломаной $v u v_1 v_2 v_3 v_4$ не превосходят $\frac{\varepsilon}{2}$, то $\lVert v - v_1 \rVert < \frac{\varepsilon}{2}$, $\lVert u - v_2 \rVert < \frac{\varepsilon}{2}$, $\lVert v - v_3 \rVert < \varepsilon$, $\lVert u - v_4 \rVert < \varepsilon$, тем самым $v_i \in Q$, $i = 1,2,3,4$.

\begin{figure}[ht]
  \centering
  \includegraphics[width=10cm]{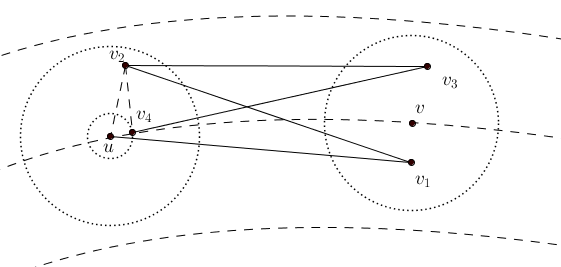}\\
  \caption{Путь длины 4 между точками $u$ и $v_4$.}
\end{figure}

\noindent При этом
$$l_1 = \lVert u - v_2 \rVert \in \left [0; 2 \sin \frac{\varepsilon}{4} \right],$$
$$l_2 = \lVert v_2 - v_4 \rVert \in \left [0; 2 \sin \frac{\varepsilon}{4} \right]$$
могут быть выбраны произвольным образом, а ориентированный угол между векторами $\overrightarrow{v_2 u}$ и $\overrightarrow{v_2 v_4}$ независимо выбирается из отрезка $[-\frac{\varepsilon}{4}; \frac{\varepsilon}{4}]$.
Зафиксируем прямую, содержащую вектор $\overrightarrow{v_2 u}$ (пусть, например, он ортогонален $u v$). Тогда всевозможные положения точки $v_4$ образуют некоторую фигуру, содежащую ромб с центром в $u$, стороной $2 \sin \frac{\varepsilon}{4}$ и углом $\frac{\varepsilon}{2}$. Следовательно, существует путь длины 4 между $u$ и любой точкой из $\gamma$-окрестности  $u$, где $\gamma =  \sin  \frac{\varepsilon}{2} \sin \frac{\varepsilon}{4}$.

\vskip+0.2cm


Таким образом, идя по кривой $\xi$ шагами размера $\gamma$ от $u$ к $v$, мы строим путь четной длины между точками $u$ и $v$, а значит, и нечетный цикл в $G(Q)$.
\qed

\subsection {Доказательство Леммы~\ref{L1}}

Покажем, что найдется круг сколь угодно малого радиуса, содержащий точки по крайней мере трех цветов.
Предположим противное: найдутся правильная раскраска плоскости и $\varepsilon>0$ такие, что все круги радиуса $\varepsilon$ содержат точки не более двух различных цветов.
Разобьем плоскость на квадраты со стороной $$\delta \leq \frac{2}{ \sqrt{10}} \varepsilon.$$
Тогда каждый такой квадрат раскрашен не более чем в два цвета.




По Утверждению~\ref{Utv2} внешняя граница любой связной двухцветной области ограничивает некоторую фигуру (иначе, соединив центры соседних граничных клеток, мы получим достаточно длинную ломаную).
Поскольку у любой связной двухцветной области внешняя граница состоит из одноцветных квадратов, диаметр области конечен.
Значит, мы можем рассмотреть такую связную двухцветную область, что ее внешняя граница ограничивает фигуру максимальной площади.
Добавим произвольный квадрат, смежный с областью снаружи и получим противоречие.

\qed

\subsection {Доказательство Следствия~\ref{Cor1}}

По Лемме~\ref{L1} для произвольного $\varepsilon > 0$ существует круг диаметра $\varepsilon$, в котором есть точки трех различных цветов.
Покажем, что какие-то точки трех различных цветов лежат на окружности радиуса не более $\varepsilon$.
Рассмотрим в круге диаметра $\varepsilon$ разноцветный треугольник $ABC$. У него есть тупой угол (иначе нам подходит описанная окружность треугольника $ABC$); не умаляя общности, это угол $A$. Рассмотрим точку $D$, такую что $\angle ADB = \angle ADC = \pi/3$.
Тогда $\angle BDC = 2\pi/3$. Заметим, что независимо от цвета точки $D$ один из треугольников  $ABD, ACD, BCD$
разноцветный. Радиус описанной окружности любого из треугольников не превосходит $$\frac{\varepsilon}{2\sin \angle D} = \frac{\varepsilon}{\sqrt 3} < \varepsilon.$$  Доказательство завершено, поскольку $\varepsilon$ можно выбрать произвольно.\qed

\subsection {Доказательство Теоремы~\ref{Th7}}

В основе доказательства лежит следующая конструкция: если существует треугольник с разноцветными вершинами $v_1,v_2,v_3$ из слойки и центром описанной окружности $u_0$, то слойка содержит окружность, которая является пересечением трех сфер радиуса 1 с центрами в $v_1,v_2,v_3$. При надлежащем выборе вершин треугольника эта окружность содержит нечетный цикл, и для ее раскраски потребуется еще три цвета, откуда имеем оценку на хроматическое число слойки, приведенную в формулировке теоремы.

Предположим существование правильной раскраски слойки $$\mathbb{R}^2\times[0, \varepsilon]^2 = \left\{(x,y,z,t) \mid x,y \in \mathbb{R},\ z,t \in [0, \varepsilon] \right\}$$ в несколько цветов. Далее потребуется вспомогательное

\begin{prop}

 Пусть $\phi (v_1,v_2,v_3)$~--- угол между двумерной плоскостью, содержащей точки $v_1$, $v_2$, $v_3$, и плоскостью $\{(0,0,z,t)\}$. Для произвольных $\varepsilon_2>0$, $\varepsilon_3>0$ в слойке найдется треугольник с вершинами $v_1$, $v_2$, $v_3$ трех различных цветов и углами $\alpha_1$, $\alpha_2$, $\alpha_3$, для которого выполнены условия:

\begin{equation}\label{ineq_phi}
 \phi(v_1,v_2,v_3) \leq \varepsilon_2;
\end{equation}

\begin{equation}\label{ineq_alpha}
\alpha_i \geq \frac{\pi}{5} - \varepsilon_3, \; i = 1,2,3.
\end{equation}

\end{prop}

\begin{proof} Выберем $\varepsilon_1<\varepsilon/2$. Пусть $M = \left\{(z_1,t_1), (z_2,t_2), \dots (z_5,t_5)\right\}$~--- вершины правильного пятиугольника, вписанного в окружность радиуса $\varepsilon_1$ с центром в точке $(\varepsilon/2,\varepsilon/2)$.  Сопоставим каждой точке $(x,y) \in \mathbb{R}^2$ вершины пятиугольника, лежащего в инфинитезимальном квадрате:
$$Q_{x, y} = \{(x,y)\} \times M.$$

Если при некоторых $x,y$ множество $Q_{x,y}$ раскрашено по крайней мере в 3 цвета, то треугольник с вершинами в разноцветных точках из  $Q_{x,y}$~--- искомый. Пусть верно обратное: для любых $x,y$ в $Q_{x,y}$ встречается не более чем 2 цвета. Тогда в раскраске $Q_{x,y}$ какой-либо цвет использован по крайней мере 3 раза; обозначим этот цвет $c(x,y)$.

\begin{figure}[ht]
  \centering
  \includegraphics[width=10cm]{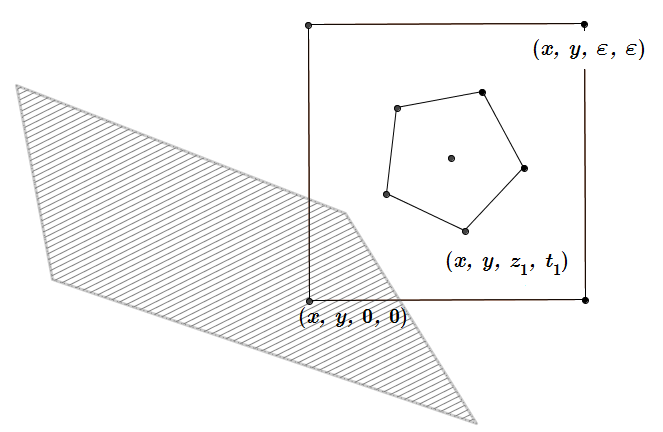}\\
  \caption{Пятерка точек $Q_{x, y} = \{(x,y)\} \times M$}
\end{figure}


Сопоставим раскраске слойки  вспомогательную раскраску плоскости  $P = \mathbb{R}^2$: точка $(x,y)$ имеет цвет $c(x,y)$. Заметим, что эта раскраска правильная. В самом деле, если точки $(x_1,y_1)$ и $(x_2,y_2)$ находятся на расстоянии 1 и раскрашены одинаково, то в соответствующих пятерках точек из слойки $Q_{x_1,y_1}=\left\{q_{1i}\right\}$, $Q_{x_2,y_2}=\left\{q_{2i}\right\}$ по крайней мере три точки имеют тот же цвет. Но поскольку $$\rho(q_{1i},q_{2i})=1, \, \,  i = 1, \dots, 5,$$
в $Q_{x_1,y_1} \cup Q_{x_2,y_2}$ может быть не более 5 точек одного цвета.

Применим Лемму~\ref{L1} к раскраске плоскости $P$: для произвольного $\delta>0$ найдутся точки $u,v,w$ трех различных цветов $c(u)$, $c(v)$, $c(w)$, попарные расстояния между которыми не превосходят $\delta$. Это означает, что в пятерках $Q_u$, $Q_v$, $Q_w$ есть по крайней мере 3 точки цветов $c(u)$, $c(v)$, $c(w)$ соответственно. Можно выбрать по одной точке из каждого набора так, чтобы они не совпадали в проекции на плоскость $(0,0,z,t)$ и имели цвета $c(u)$, $c(v)$, $c(w)$. Нетрудно проверить, что если выполнены условия
$$16\left(\frac{\delta}{\varepsilon_1}+2 \frac{\delta^2}{\varepsilon_1^2}\right) \leq \sin \varepsilon_2; $$
$$\delta \leq \frac{\varepsilon_1}{2} \sin \frac{\varepsilon_3}{2},$$
то справедливы неравенства (\ref{ineq_phi}), (\ref{ineq_alpha}).
\end{proof}

Теперь мы готовы начать доказательство Теоремы~\ref{Th7}. Рассмотрим точки $v_1, v_2, v_3$, удовлетворяющие условиям Утверждения 4. Пусть $u_0$~--- центр описанной окружности треугольника $v_1 v_2 v_3$; $\overline{n}$~--- некоторый единичный вектор, ортогональный двумерной плоскости, в которой лежит треугольник; $u_1 = u_0 + \delta_1 \overline{n}$; $L(u_1, v_1, v_2, v_3)$ --- гиперплоскость, проходящая через точки $u_1$, $v_1$, $v_2$, $v_3$.
Пусть $B(u_1; \delta_2) \subset L(u_1, v_1, v_2, v_3)$~--- трехмерный открытый шар радиуса $\delta_2 > 0$ с центром в $u_1$.



Для точки $w \in B(u_1; \delta_2)$ определим
$$T_{1}(w) = S(v_2; 1) \cap S(v_3; 1) \cap S(w; 1),$$
где $S(v;1)$~--- сфера единичного радиуса с центром $v$.

\begin{figure}[ht]
  \centering
  \includegraphics[width=12cm]{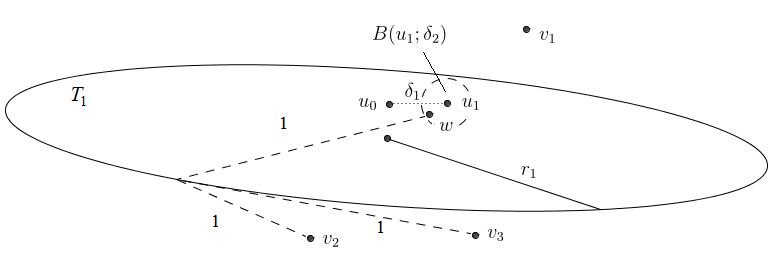}\\
  \caption{Построение окружности запрещенного радиуса}
\end{figure}


Пусть радиус окружности $T_{1}(w)$ равен $r_1(w)$; окружности $T_{2}(w)$, $T_{3}(w)$ и их радиусы $r_2(w),r_3(w)$ определены аналогично.

Поскольку вершины треугольников $w v_1 v_2$, $w v_2 v_3$, $w v_1 v_3$ находятся на расстоянии не более $\delta +\delta_1 +\delta_2$ от вершин треугольника, лежащего в сечении $\{(0,0)\} \times [0, \varepsilon]^2$, то при достаточно малых $\delta$, $\delta_1$, $\delta_2$ соответствующие окружности лежат внутри слойки.

 Для определенного выше трехмерного шара $B(u_1; \delta_2)$ введем функцию
$$r: B(u_1; \delta_2) \rightarrow \mathbb{R}^3;$$
$$r(w) = (r_1(w),r_2(w),r_3(w)).$$

Заметим, что при $w = u_1$ градиенты функций $r_i(w)$ коллинеарны медианам равнобедренных треугольников $u_1 v_2 v_3$, $u_1 v_1 v_3$, $u_1 v_1 v_2$, проведенным из вершины $u_1$:
$$
\nabla r_1(u_1) = \lambda_1 \left(u_1 - (v_2+v_3)/2 \right); \quad
\nabla r_2(u_1) = \lambda_2 \left(u_1 - (v_1+v_3)/2 \right); \quad
\nabla r_3(u_1) = \lambda_3 \left(u_1 - (v_1+v_2)/2 \right),
$$
причем $\lambda_i \neq 0$ и симплекс $u_1 v_1 v_2 v_3$ невырожден. Следовательно, при $w = u_1$ якобиан $\partial r / \partial w$ отличен от 0, и в окрестности $u_1$ для $r(\cdot)$ выполнены условия теоремы об обратной функции. Но поскольку запрещенные радиусы всюду плотны в окрестности каждого из значений $r_1(w)$, $r_2(w)$, $r_3(w)$, то найдется тройка запрещенных радиусов $r^*_1$, $r^*_2$, $r^*_3$, имеющая прообраз $u^*$ в $B(u_1; \delta_2)$.

 Тогда при любом цвете точки $u^*$ хотя бы один из треугольников $u^*v_1v_2$, $u^*v_2v_3$, $u^*v_1v_3$ будет иметь разноцветные вершины, а на соответствующей окружности запрещенного радиуса должны встречаться еще по крайней мере 3 цвета. Тем самым $\chi(\mathbb{R}^2\times[0, \varepsilon]^2) \geq 6$. \qed

\subsection {Доказательство Теоремы \ref{Th8}}

Идея доказательства заключается в построении семейства замкнутых множеств, диаметр каждого из которых не превосходит 2, и которые также образуют покрытие $\mathbb{R}^n$. После того, как искомое семейство найдено, утверждение следует непосредственно из определения топологической размерности. Используется стандартная топология $\mathbb{R}^n$.

Напомним, что $C_i$~--- множество всех точек $\mathbb{R}^n$, раскрашенных в $i$-й цвет, $1\leq i \leq m$, и положим
$${C^*_i:=\overline{\Int \, \overline{C_i}}}\qquad \mbox{(замыкание\ внутренности\ замыкания)}.$$

\noindent Разобьем каждое из множеств $C^*_i$ на компоненты связности (в смысле стандартной топологии): $$C^*_i = \bigcup_{\alpha\in A_i} D_{\alpha}.$$

\noindent Для краткости введем обозначение $\{ D_\alpha \} = \bigcup\limits_{i=1}^m \bigcup\limits_{\alpha \in A_i} D_\alpha$.

\medskip

\noindent{\bf (i).} {\it Множества $C^*_i$  образуют покрытие $\mathbb{R}^n$.}

\medskip

\noindent Пусть верно обратное: $\exists v: \forall i \ \ v\notin C^*_i$. Тогда существует открытый шар $B(v; \varepsilon)$:
 $$B(v; \varepsilon)\cap C^*_i =\emptyset; \ \ B(v; \varepsilon)\subset\bigcup C_i.$$
 Рассмотрим некоторый шар $$B^1 \subset B(v; \varepsilon) \setminus \overline{C}_1.$$ Очевидно, $B^1$ не может быть подмножеством $\overline{C}_i$~--- иначе пересечение внутренности $\overline{C}_i$ и $B(v; \varepsilon)$ было бы непусто. Определим последовательность вложенных шаров $$B^{k+1} \subset B^k \setminus \overline{C}_k.$$ Точки из $B^{m+1}$  не принадлежат ни одному из $\overline{C}_i$, что противоречит исходным предположениям.

\medskip

\noindent{\bf (ii).} {\it Если сфера $S$ радиуса 1 с центром в точке $v$ содержит внутренние точки множеств $\overline{C}_i$, ${1\leq i\leq k \leq n}$, то $v$ принадлежит  хотя бы одному из множеств $C^*_j$, $k+1\leq j\leq m$.}

\medskip

\noindent

Можно выбрать точки $x_1, ... , x_n$ таким образом, что
$$x_i\in S  \cap \Int \overline{C}_i, \quad \ 1\leq i \leq k;$$
$$ x_i \in S, \quad k+1 \leq i \leq n;$$
и $\{v, x_1, ... , x_n\}$ находятся в общем положении (в смысле невырожденности симплексов). Определим такое $\varepsilon > 0$, что $B(x_i; \varepsilon) \subset \overline{C}_i$, $1 \leq i \leq k$.
Цвет точки
$$w = w(q_1,\dots,q_k) = \bigcup\limits_{1 \leq i \leq k} S(q_i; 1),$$
 если она определена, отличается от каждого из цветов точек $q_1, ... , q_n$.
Пусть $$z \in B(0; \varepsilon); \quad y_i = x_i + z.$$  В достаточно малой окрестности набора точек $y_i$ функция $w(\cdot)$ определена и непрерывна по каждому из аргументов. Выберем такие точки
$$y'_i \in C_i, \quad 1 \leq i \leq k;$$
$$y'_i = y_i, \quad 1 \leq k+1 \leq n,$$ что существует $w(y'_1, ... , y'_n)$. Тогда
$$w(y'_1,\dots,y'_k) \in \bigcup_{j=k+1}^{m} {C}_j.$$
При этом
$$\delta(y'_1, ... ,y'_k) = \max \limits_{1 \leq i \leq k} \lVert y'_i - y_i \rVert$$
может быть сколь угодно малым, следовательно
$$w(y_1,\dots,y_k) \in \bigcup_{j=k+1}^{m} \overline{C}_j.$$

\noindent Поскольку $z \in B(0; \varepsilon)$ произвольно,
$$B(v; \varepsilon) \subset \bigcup_{j=k+1}^{m} \overline{C}_j.$$
Следовательно, хотя бы одно из  множеств $\overline{C}_j$, $j = k+1, ... ,m$ всюду плотно в некоторой окрестности $v$.

\medskip

\noindent{\bf (iii).}\ {\it Если некоторая точка $v\in\mathbb{R}^n$ покрыта не более чем $n$ множествами из $\{D_\alpha\}$, то диаметр хотя бы одного из этих множеств не превосходит $2$.}

\medskip

\noindent В противном случае каждое множество из $\{D_\alpha\}$, покрывающее $v$, имеет непустое пересечение со сферой $S$ радиуса 1 с центром в точке $v$. Без ограничения общности предположим, что точку $v$ покрывают множества $D_1,\dots,D_n$, являющиеся компонентами связности $C^*_1,\dots,C^*_n$ соответственно. Пусть $\operatorname{min}\{\Diam D_i\} = 2+\delta$.



Пусть $w \in \mathbb{R}^n, \, \|w\| = 1$~--- некоторое направление; $S_\eta$~--- сфера радиуса $1$ с центром в точке $u(\eta) = v+\eta w$, $\eta\in\mathbb{R}_+$.
Тогда множество $$T_i = \{\eta\in\mathbb{R}_+: \ \ S_\eta \cap\Int D_i \neq \emptyset; \ 1 \leq i \leq n \}$$ всюду плотно на отрезке $[0,\delta]$.
Следовательно,  $$\forall \eta \in [0,\delta] \ \ u(\eta) \in \bigcup^m_{n+1} \overline{C}_j\ .$$
Аналогичное рассуждение справедливо для произвольного  вектора $w$ единичной длины.
Но тогда любая окрестность $v$ содержит шар, являющийся подмножеством $\bigcup^m_{n+1} \overline{C}_j$, а следовательно, и внутреннюю точку хотя бы одного из множеств $\overline{C}_j, \ j>n$.  Полученное противоречие доказывает (iii).

\medskip

\noindent{\bf (iv).}\ {\it Если любая точка $\mathbb{R}^n$ покрыта не более чем $n$ множествами из $\{D_\alpha\}$, то семейство множеств $\Delta = \{D_\alpha | \ \Diam(D_\alpha)\leq 2\}$ покрывает $\mathbb{R}^n$}.

\medskip

\noindent Очевидным образом следует из (iii).

\medskip

\noindent{\bf (v).}\ {\it Найдутся множества $D'_1, D'_2,\dots,D'_{n+1} \in \Delta$, имеющие непустое пересечение.}

\medskip

\noindent Рассмотрим шар $B(0; R)\subset\mathbb{R}^n$ и его покрытие  семейством множеств $\Delta$. Согласно одному из  определений топологической размерности (см. П.~С.~Александров и Б.~А.~Пасынков~\cite{AlPas}) при достаточно большом $R$ среди покрывающих $B(0; R)$ замкнутых множеств, диаметр каждого из которых не превосходит $2$, найдутся $n+1$, пересечение которых непусто.

\medskip

\noindent{\bf (vi).}\ {\it Найдется $n+1$ множество из семейства $\{C_i\}$, пересечение замыканий которых непусто.}

\medskip

\noindent Пусть для множеств $D'_1, D'_2,\dots,D'_{n+1} \in \Delta$ выполнено
$$
\bigcap \limits_{i=1}^{n+1} D'_i \neq \emptyset;
$$
$$
D'_i \subset C^*_{l_i}, \, \, i=1, 2, ... , n+1.
$$
Заметим, что индексы $l_i$ попарно различны, поскольку в противном случае попарно пересекающиеся множества $D'_i$ не могли бы являться различными компонентами связности $C^*_i$. Следовательно,
$$
\emptyset \neq \bigcap \limits_{i=1}^{n+1} D'_i \subset  \bigcap \limits_{i=1}^{n+1} C^*_{l_i} \subset \bigcap \limits_{i=1}^{n+1} \overline{C}_{l_i},
$$
и $\{C_{l_i}\}$~--- искомое подсемейство.  \qed

\medskip
\textbf{Замечание.} Пользуясь леммой Шпернера, нетрудно получить оценку на радиус шара, который содержит хотя бы одну точку, принадлежащую $n+1$ из множеств $\overline{C}_i$.






\subsection {Доказательство Теоремы \ref{Th9}}

Пусть координата $x$ полноценная, а координаты $y$, $z$ и $t$ инфинитезимальные.

\noindent{\bf Верхняя оценка.}\ Покрасим точки с координатами $(x,y,z,t)$ при $\frac{2k}{3} < x \leq \frac{2(k+1)}{3}$ в цвет $k \ \mbox{mod}\ 3$ ($k$~--- целое число).

\noindent{\bf Нижняя оценка.}\ Покажем наличие нечетного цикла в дистанционном графе $G_1 (\mathbb{Q}\times [0,\varepsilon]_\mathbb{Q}^3)$.

Рассмотрим четное $n > 2\varepsilon^{-2}$ и вектор $e = (1-n^{-1},bn^{-1},cn^{-1},dn^{-1})$, такой что $b^2+c^2+d^2=2n-1$. Этот вектор имеет единичную длину.
Заметим, что $e$ умещается в нашу полосу, поскольку
$$\max (|b|n^{-1}, |c|n^{-1}, |d|n^{-1}) < \sqrt{\frac{2}{n}} < \varepsilon.$$
Также рассмотрим вектор $e'=(1-n^{-1},-bn^{-1},-cn^{-1},-dn^{-1})$ и последовательность точек $A_i$, определенную следующим образом
$$A_0 := (0,0,0,0); \ \ \ A_{2k+1} := A_{2k}+e; \ \ \ A_{2k+2} := A_{2k+1}+e'.$$
Нетрудно видеть, что $A_n=(n-1,0,0,0)$, поскольку $n$ четно. Таким образом, точки $A_0,\dots,A_n$ и точки $(1,0,0,0),\dots,(n-2,0,0,0)$ образуют искомый нечетный цикл.
Остается заметить, что для любого $\varepsilon>0$ найдутся целые $n$, $b$, $c$, $d$, удовлетворяющие нашим условиям (например, $b = c = d = 2l+1$, $n = 6l^2 + 6l + 2$ при достаточно большом $l$). \qed




\section {Заключение и вопросы для дальнейшего исследования}

Мы показали, что $5 \leq \chi(\mathbb{R}^2 \times [0,\varepsilon]) \leq 7$ и ${6 \leq \chi(\mathbb{R}^2 \times [0,\varepsilon]^2) \leq 7}$, а также $\chi(\mathbb{R}^2 \times [0,\varepsilon]^k) \leq 7$ при достаточно маленьком $\varepsilon$. Естественным образом возникает вопрос о существовании такого $k$, что \\
$\chi(\mathbb{R}^2 \times [0,\varepsilon]^k) = 7$ при произвольно малом $\varepsilon$.

Как можно заметить, в тех случаях, когда мы можем посчитать хроматическое число вещественной слойки, имеет место дискретная непрерывность.
Выполняется ли это свойство в общем случае? Иными словами, является ли функция $\chi(\mathbb{K}^n \times [0,\varepsilon]^m)$, где $\mathbb{K} \in \{\mathbb{R}, \mathbb{Q}\}$ дискретно непрерывной по $\varepsilon$?

В одномерном и двумерном случае добавление одной инфинитезимальной размерности увеличивает нижнюю оценку на хроматическое число пространства.
Из общих соображений следует, что должна выполняться оценка $\chi(\mathbb{R}^3\times [0,\varepsilon]) \geq 7$ , где $\varepsilon$~--- произвольное положительное число, однако нам не удалось это доказать. Более того, хочется предположить, что $\chi (\mathbb{R}^n \times [0, \varepsilon]) > \chi (\mathbb{R}^n)$, однако это уже куда более сложное утверждение.

А.~Б.~Купавский в работе~\cite{K} поставил вопрос о максимальном гарантированном количестве цветов на $m$-мерной сферы радиуса $r$ при правильной раскраске $n$-мерного
пространства в конечное количество цветов. В той же работе были получены оценки при $r$, отделенных от нуля. Лемма~\ref{L1} дополняет эти результаты при инфинитезимальных $r$, но только в случае $n = 2$, $m = 1$. По-видимому, из Теоремы~\ref{Th8} можно вывести аналогичный результат для $n = m+1$ при произвольном $n>2$.

Также интерес представляет и обратная задача.
По натуральному числу $k$  требуется построить ``разумное'' пространство с хроматическим числом ровно $k$.
Например, было бы очень интересным такое утверждение для пространства с большим аффинным подпроcтранством, в частности для $[0,h_1] \times\dots\times [0,h_m] \times [0,\varepsilon]^l \times \mathbb{R}^s$ при $s>0$.

\vskip 0.5cm

\textbf{Благодарности.}
Работа была поддержана Российским Научным Фондом:
Теоремы 7 и 11, а также Лемма 1 грантом 16-11-10039, а Теоремы 8 и 9 "--- грантом 17-11-01377.
Авторы благодарят Мишу Баска и Федора Петрова за плодотворное обсуждение и ценные замечания, а также анонимного рецензента за указание огромного количества неточностей и предложения по их исправлению.

\vskip 0.5cm

\textbf{Сведения об авторах.}

Канель-Белов Алексей Яковлевич, Лаборатория продвинутой комбинаторики, Московский физико-технический институт, Институтский пер., 9, Долгопрудный, Московская обл., 141701.

kanelster@gmail.com

Воронов Всеволод Александрович, Институт динамики систем и теории управления
имени В.М. Матросова, Сибирского отделения РАН,
664033, Иркутск, ул. Лермонтова, 134, а/я 292.

v-vor@yandex.ru

Черкашин Данила Дмитриевич, Лаборатория продвинутой комбинаторики, Московский физико-технический институт, Институтский пер., 9, Долгопрудный, Московская обл., 141701; Лаборатория им. П.Л. Чебышева, Санкт-Петербургский государственный университет, 14 линия В.О., дом 29Б, Санкт-Петербург 199178 Россия;
С.-Петербургское отделение Математического института
им. В. А. Стеклова РАН
191023, Санкт-Петербург
наб. р. Фонтанки, 27
Россия

matelk@mail.ru



\end{document}